\documentclass[12pt,twoside]{amsart}
  \usepackage{amssymb}
  \usepackage{amscd}

\newcommand{\1}{^{-1}}
\newcommand{\tensor}{\otimes}
\newcommand{\into}{\hookrightarrow}
\DeclareMathSymbol{\subsetneq}{\mathrel}{AMSb}{"28}

\newcommand{\A}{\mathbb A}
\newcommand{\C}{\mathbb C}
\newcommand{\G}{\mathbb G}
\newcommand{\proj}{\mathbb P}
\newcommand{\sS}{\mathcal S}
\newcommand{\Spec}{\operatorname{Spec}}
\newcommand{\codim}{\operatorname{codim}}
\newcommand{\wave}{\widetilde}
\newcommand{\Oh}{\mathcal O}

\newcommand{\ep}{\varepsilon}
\newcommand{\si}{\sigma}
\newcommand{\Si}{\Sigma}

\theoremstyle{plain}
\newtheorem{thm}{Theorem}
\newtheorem{lemma}{Lemma}
\theoremstyle{definition}
\newtheorem{claim}{Claim}

\title[Bogomolov--Pantev resolution]{The Bogomolov--Pantev resolution,
  an expository account}
\author{Kapil H. Paranjape}
\thanks{The author thanks the University of Warwick for hospitality
during the period when this work was done.}
\address{Institute of Mathematical Sciences, CIT Campus, Tharamani,
Chennai 600 113, INDIA.}
\email{kapil@imsc.ernet.in}

\begin{document}
 \maketitle

\section*{Statement of the result}

Bogomolov and Pantev~\cite{Bog-Pant} have recently discovered a rather
elegant geometric proof of the weak Hironaka theorem on resolution of
singularities:
 \begin{thm} \label{main}
   Let $X$ be a projective variety and $Z$ a proper Zariski closed
   subscheme of $X$. There is a projective birational map
   $\ep\colon\wave{X}\to X$ such that $\wave X$ is smooth, the scheme
   theoretic inverse image $\ep\1(Z)$ is a Cartier divisor, the
   exceptional locus of $\ep$ is a divisor and the union of these two
   divisors is a divisor with simple normal crossings.
 \end{thm}
Before the work cited above, and that of Abramovich and de Jong
\cite{Ab-Jong} (appearing at roughly the same time) the only proof of
this theorem was as a corollary of the famous result of Hironaka
\cite{Hironaka}. These new proofs were inspired by the recent work of
de Jong \cite{Jong}, which Bogomolov and Pantev combine with a
beautiful idea of Belyi \cite{Belyi} ``simplifying'' the ramification
locus of a covering of $\proj^1$ by successively folding up the
$\proj^1$ onto itself, over a fixed base. This latter step
unfortunately only works in characteristic zero, limiting the scope
of the argument (Abramovich and de Jong's paper gives some results
even in characteristic $p$). Hence, we work over the field of complex
numbers; the argument also works (with suitable modifications about
rationality) over any field of characteristic zero.

The outline of the argument we follow is the same as that of the paper
of Bogomolov and Pantev; however we offer different (and we hope
simpler) proofs of the corresponding lemmas. To begin with, their
argument using Grassmannians is replaced by an application of Noether
normalisation in Section~\ref{bundle}. Belyi's argument to reduce the
degree of individual components of the ramification locus is presented in
purely algebraic form in Section~\ref{bundle}, Lemmas~\ref{belyi1}
and~\ref{belyi2}. The presentation of Bogomolov and Pantev refers to
``semi-stable families of pointed curves of genus 0''; the precise
result required from that theory is proved here by means of blowups in
Section~\ref{fibrations}. Finally we also give a summary (in
Section~\ref{toric}) of the desingularisation of toroidal embeddings
which is used by both the papers \cite{Bog-Pant,Ab-Jong}. To summarise,
the aim of this account is to give all the details of the argument, so
that it should be accessible to anyone with a basic knowledge of
algebraic geometry (as contained for example in Mumford's book
\cite{Redbook}).

\subsection*{Acknowledgements} The last statement above was the
``challenge'' posed by Miles Ried who has contributed to this paper in
numerous ways. Discussions with Lawrence Ein and D.~S.~Nagaraj have
greatly helped me in understanding this proof better. V.~S.~Sunder
patiently heard out my ramblings on ``convex'' barycentric
subdivisions. Madhav Nori pointed out a flaw in an earlier argument in
the section on genus 0 fibrations. Rob Lazarsfeld suggested that the
statement of the main result be strengthened---this did not require a
change in the proof! I thank the referee for reading the paper very
carefully and providing a number of suggestions. Of course, many
thanks to Fedor Bogomolov and Tony Pantev for their beautiful proof!

\section*{Some reductions}
For completeness, we define the notion of {\em exceptional locus} that
we use. For a proper birational morphism $\ep:Y\to X$ of varieties,
the {\em exceptional locus} is defined as the complement of the
largest $\ep$-saturated open subset of $Y$ on which $\ep$ is an
isomorphism. We note that the exceptional locus of the composite
$f\circ g$ is the union of the inverse image under $g$ of the
exceptional locus for $f$ and the exceptional locus for $g$.

Let $f:X_1\to X$ be the blowup of $X$ along the scheme $Z$; we thus
ensure that $f^{-1}(Z)$ supports a Cartier divisor. Now let $g:X_2\to
X_1$ be the blowup along the reduced scheme supported on the
exceptional locus $E(f)$ of $f$. Finally, let $h:X_3\to X_2$ be the
blowup along the reduced scheme supported on the singular locus of
$X_2$. Let $Z'$ denote the union of the supports of inverse images of
$Z$, $E(f)$ and $(X_2)_{\rm sing}$ in $X_3$.
\begin{lemma}
  With notation as above, let $\ep_1:\tilde{X}\to X_3$ be a projective
  birational morphism with $\tilde{X}$ smooth such that $\ep_1\1(Z')$
  is contained in a divisor with simple normal crossings. The composite
  morphism $\ep=f\circ g\circ h\circ\ep_1$ then satisfies the conclusion
  of Theorem~\ref{main}.
\end{lemma}
\begin{proof}
  Since any sub-divisor of a divisor with simple normal crossings (see
  Section~\ref{fibrations} for the definition) is again a divisor with
  simple normal crossings it is enough to show that each of $\ep^{-1}(Z)$
  and the exceptional locus $E(\ep)$ are divisors.

  First of all, $f^{-1}(Z)$ supports a Cartier divisor, thus so does
  $\ep^{-1}(Z)$. It follows that the latter is a divisor.
  
  If $D$ is any divisorial component of $(X_2)_{\rm sing}$, then $X_2$
  is not normal along $D$. It follows that $X_3\to X_2$ is {\em not}
  an isomorphism along $D$ and thus $h^{-1}(D)$ is contained in the
  exceptional locus for $h$. The inverse image of any codimension 2
  component of $(X_2)_{\rm sing}$ is clearly contained in the
  exceptional locus of $h$. Thus, the exceptional locus for $h$ is
  precisely $h^{-1}((X_2)_{\rm sing})$. The exceptional locus of
  $f\circ g:X_2\to X$ is $g^{-1}(E(f))$ which supports a Cartier
  divisor by construction of $g$. It follows that the exceptional
  locus of $f\circ g\circ h:X_3\to X$ supports a Cartier divisor.
  Moreover, the complement $U$ of this divisor in $X_3$ is smooth and
  is mapped isomorphically to an open subset of $X$ by $f\circ g\circ
  h$. Let $V=\ep_1^{-1}(U)$. The morphism $\ep_1:V\to U$ is a proper
  birational morphism of smooth varieties and thus its exceptional
  locus is a divisor $D$ in $V$. It thus follows that the exceptional
  locus of $\ep$ is a divisor.
\end{proof}
As a corollary we see that the main result Theorem~\ref{main} is
equivalent to the apparently weaker result:
\begin{thm}\label{weaker}
  Let $X$ be a projective variety and $Z$ a proper Zariski closed
  subset of $X$. There is a projective birational map
  $\ep\colon\wave{X}\to X$ such that $\wave X$ is smooth. Moreover,
  the set theoretic inverse image $\ep\1(Z)$ and the exceptional locus
  of $\ep$ are contained in a divisor with simple normal crossings.
\end{thm}
We will now prove the two equivalent theorems, by proving
Theorem~\ref{weaker}. We can make the following trivial reductions:
 \begin{enumerate}
 \item We may assume that $X$ is normal by replacing $X$ with its
  normalisation and replacing $Z$ by the union of its inverse image
  with the exceptional locus of the normalisation.
\item Any normal curve is smooth and any proper closed subset of a
  smooth curve is a simple normal crossing divisor~(!). Hence the
  result is true in dimension 1.
\item We prove the theorem by induction on the dimension of $X$. Thus
  we may assume that $n=\dim X>1$ and that Theorem~\ref{main} 
  holds for all varieties of smaller dimension than $X$.
\item In order to prove Theorem~\ref{weaker} for a given $Z$ it is enough to
  prove it for a larger $Z'\supset Z$, as long as $Z'$ is a proper
  subset of $X$. Thus, we can expand $Z$ at any stage in our
  construction of the resolution. In particular, whenever we blow up
  $X$ we will expand $Z$ to include the exceptional locus of the blow
  up. 
 \end{enumerate}

\section{$\proj^1$-bundles}\label{bundle}
We will prove the following:
 \begin{claim}\label{claim1}
   After replacing $X$ by $X'$, where the latter is the blow up of $X$
   at a finite set of smooth points, and replacing $Z$ by the union of
   its inverse image in $X'$ and the exceptional divisor of the blow
   up, we have:
 \begin{enumerate}
 \item There is a finite surjective map $f\colon X\to\proj_Y(F)$,
   where $Y$ is smooth and $F$ is a rank $2$ vector bundle on $Y$.
 \item The morphism $f$ is finite and \'etale outside a divisor
   $B=\bigcup_{i=1}^n s_i(Y)$ which is the union of a finite number of
   sections $s_i\colon Y\to\proj_Y(F)$ of the $\proj^1$-bundle
   $\proj_Y(F)$.
 \item The image $f(Z)$ of $Z$ under $f$ contained in $B$.
 \end{enumerate}
 \end{claim}

 \begin{lemma}\label{noether}
After replacing $X$ by its blowup at a suitable finite set of smooth
points, we have the following situation. There is a smooth variety $Y$, a
line bundle $L$ and a finite surjective map $g\colon X\to Q=
\proj_Y(\Oh\oplus L)$ such that if $E$ denotes the section corresponding
to the quotient homomorphism $\Oh\oplus L\to\Oh$, then there is a divisor
$B\subset Q$ disjoint from $E$ and containing $g(Z)$ and the branch locus
of $g$. Moreover, $E$ is the image of the exceptional locus of the
blowup of $X$ at the finite set of smooth points.
 \end{lemma}

 \begin{proof}
This is essentially just Noether normalisation. Let $X\into\proj^N$ be a
projective embedding of $X$ and $L=\proj^{N-n-1}$ a linear subspace in
$\proj^N$ that does not meet $X$, where $n=\dim X$. The projection from
$L$ gives a finite map $X\to\proj^n$. Let $B$ be a hypersurface in
$\proj^n$ that contains the image of $Z$ and the branch locus of
the map. Let $p$ be a point of $\proj^n$ not on $B$. Replace $X$ by its
blowup at the points lying over $p$ (which are all smooth) and let
$Y=\proj^{n-1}$ with $L=\Oh(1)$. The blowup of $\proj^n$ at $p$ is
naturally isomorphic to $\proj_Y(\Oh\oplus L)$ and the resulting
morphism from $X$ to this $\proj^1$-bundle has the required form. 
 \end{proof}

Now $B$ is a divisor in the $\proj^1$-bundle $Q\to Y$ that does not meet
a section $E$. Thus the projection $B\to Y$ is finite and flat. In fact
we have:

\begin{lemma}\label{belyi1}
In the above situation, if $\Oh_Q(1)$ denotes the universal quotient
bundle, then $\Oh_Q(B)=\Oh_Q(b)$ where $b$ is the degree of the map
$B\to Y$.
\end{lemma}

\begin{proof}
Consider the line bundle $\Oh_Q(B)\otimes\Oh_Q(-b)$. Since this line
bundle is trivial on the fibres of $Q\to Y$, it is the pullback of a line
bundle from $Y$. But it restricts to the trivial bundle on $E$, which is
a section of $Q\to Y$. Hence it is the trivial line bundle.
\end{proof}
Set $Q=\proj_Y(\Oh\oplus L)$, and let $B$ be a divisor in $Q$ which is
finite over $Y$; assume moreover that $B=B_1\cup\bigcup_{i=1}^r s_i(Y)$,
where $s_i$ are sections, and $B_1$ does not meet $E$. Now let $d(B)$
denote the maximum degree over $Y$ of any irreducible component of $B_1$
and $m(B)$ the number of components of this degree. We wish to construct
a new map $X\to Q'=\proj_Y(\Oh\oplus L^N)$ for which $B_1$ is empty.
Thus we may assume that $d(B)>1$ and $m(B)>0$. The following lemma
shows that we can arrange for at least one of these numbers to drop, and
that completes the required inductive step.

\begin{lemma}\label{belyi2}
There is a map $h\colon Q\to Q'=\proj_Y(\Oh\oplus L^d)$ such that if $B'$
is the union of $h(B)$ and the branch locus of $h$ then
$(d(B),m(B))>(d(B'),m(B'))$ in the lexicographic ordering. Moreover,
under this map the image of the divisor $E$ in $Q$ is the
corresponding divisor $E'$ in $Q'$.
 \end{lemma}

\begin{proof}
Let $A$ be an irreducible component of $B_1$ of degree $d=d(B)$ and
consider the two maps $\Oh_Q\to\Oh_Q(d)$ and $L_Q^d\to\Oh_Q(d)$ given by
$A$ and $d\cdot E$. Since these two divisors do not meet, the direct sum
$\Oh_Q\oplus L_Q^d\to\Oh_Q(d)$ is surjective. Thus, by the universal
property of $Q'$ we obtain a morphism $Q\to Q'$ which on every fibre is
a map $\proj^1\to\proj^1$ of degree $d$. Its ramification divisor $R$
has the form $R=(d-1)E + R'$ for some divisor $R'$ that does not meet
$E$. Moreover, the Hurwitz formula on the general fibre $\proj^1$
(equivalently, computing the canonical divisors) gives
$\Oh_Q(R')=\Oh_Q(d-1)$. Thus we have the required result.
\end{proof}
We compose the morphism obtained from Lemma~\ref{noether} with a
succession of morphisms obtained from Lemma~\ref{belyi2}. By the
latter lemma the pair $(d(B),m(B))$ can be reduced until eventually
$B_1$ becomes empty and thus the composite morphism $f$ is as stated
in Claim~\ref{claim1}.

\section{Genus 0 fibrations}\label{fibrations}

We will prove the following:
\begin{claim}\label{claim2}
  We can replace $X$ by a blow up $X'$ and $Z$ by the union of its
  inverse image in $X'$ with the exceptional locus of the blow up, so
  that:
 \begin{enumerate}
 \item There is a finite map $f\colon X\to W$ with $W$ smooth.
 \item The union of the image $f(Z)$ of $Z$, the exceptional locus $E$
   and the branch locus of $f$ is contained in a divisor $D$ with
   simple normal crossings (or strict normal crossings).
 \end{enumerate}
 \end{claim}
For completeness, we recall the definition of a simple normal crossing
divisor $D$ in a smooth projective variety. If $D=\bigcup_{i=1}^n D_i$
is the decomposition of the divisor into irreducible components, then
for each $I\subset\{1,\ldots,n\}$ the scheme theoretic intersection
$D_I=\bigcap_{i\in I} D_i$ is reduced and smooth of codimension equal
to the cardinality $\#I$ of the set $I$.

Write $g\colon X\to P=\proj_Y(F)$ for the map constructed in
Section~\ref{bundle} and let $\{s_i\}_{i=1}^m$ be sections of $P\to Y$
such that the union $\bigcup_{i=1}^m s_i(Y)$ contains $g(Z)$, the
branch locus of $g$ and the image of the exceptional locus of the
blowups already performed. Let
 \[
 Z_1 = \bigcup_{1\le i < j\le m} p_Y(s_i(Y)\cap s_j(Y))
 \]
be the divisor in $Y$ obtained as the image of the pairwise
intersections of the sections. We apply the induction hypothesis 0.3
to obtain a map $\ep_Y\colon\wave Y\to Y$ such that the inverse image
$\ep_Y\1(Z)$ is a divisor with simple normal crossings in $\wave
Y$. Moreover, note that $Y$ as constructed in the previous section is
$\proj^{n-1}$ and thus normal. Thus, we can also assume (by induction)
that the exceptional locus $E_Y$ of $\ep_Y$ is a divisor, such that
the union $Z_1=E_Y\cup\ep_Y\1(Z)$ is also a divisor with simple normal
crossings. 

Replacing $Y$ by $\wave Y$ and $P$ by its pullback, we obtain a
configuration with the following properties (A):

 \begin{enumerate}
 \renewcommand{\labelenumi}{(\arabic{enumi})}

\item $p\colon P\to Y$ is a flat morphism of smooth varieties whose
reduced fibres are trees of projective lines $\proj^1$ (in other words,
$p$ is a {\em genus $0$ fibration});

\item we have a finite collection $\{s_i\}_{i=1}^m$ of sections of $p$
such that
 \[
 \bigcup_{1\le i < j\le m} p_Y(s_i(Y) \cap s_j(Y)) \subset Z_1
 \]
where $Z_1$ is a divisor with simple normal crossings in $Y$;

\item the morphism $p$ is smooth outside $p\1(Z_1)$ and the latter
 is a divisor with simple normal crossings.
\end{enumerate}
We wish to perform a succession of blowups resulting in a configuration
still having the same properties, but with the sections disjoint.

For each component $C$ of $Z_1$ and each section $s_i$, the image
$C'=s_i(C)\subset P$ is a codimension $2$ subvariety. Write $n(C')$ for
the number of $j$ such that $s_j(C)=C'$, and $n_P$ for the maximum of
such $n(C')$.

\begin{lemma}
 Let $C'$ be so chosen that $n(C')=n_P$. Then for any $j$, either
$s_j(C)=C'$ or $s_j(Y)\cap C'=\emptyset$.
\end{lemma}

\begin{proof}
Suppose $j$ is such that $s_j(Y)\cap C'$ is a proper nonempty subset of
$C'$, and let $C''$ be an irreducible component of this intersection.
Choose $i$ so that $s_i(C)=C'$, let $D'$ be a component of $s_i(Y)\cap
s_j(Y)$ containing $C''$, and denote by $D$ the image $p(D')$. On the one
hand we have
 \[
 s_j(C)\cap C' \subset s_j(Y)\cap C' \subsetneq C' = s_i(C),
 \]
while $D'=s_i(D)=s_j(D)$. Thus we see that $D$ and $C$ are different
components of $Z_1$ containing $p(C'')$. But $C''\subset C'$ has
codimension 3 in $P$, and $p$ restricts to an isomorphism from $C'$ to
$C$. Thus $p(C'')$ is of codimension 2 in $Y$; hence it is an irreducible
component of $D\cap C$. Since $Z_1$ is a simple normal crossing divisor
we see that $D$ and $C$ are the {\em only} components of $Z_1$ that
contain $p(C'')$.

Take any $k$ such that $s_k(C)=C'$; then $s_k(Y)\cap s_j(Y)$ contains a
component $E'$ which contains $C''$. Then by the above reasoning we must
have $D=p(E')$, and thus
 \[
 D' = s_j(D) = s_j(p(E')) = E' = s_k(p(E')) = s_k(D)
 \]
But then we get $s_k(D)=D'$ for all $k$ such that $s_k(C)=C'$, and in
addition, $s_j(D)=D'$ so that $n(D')=n(C')+1$. This contradicts the
maximality of $n(C')$.

\end{proof}

Together with the preceding lemma, the following result shows that
blowing up $C'$ with $n(C')=n_P$ again leads to a configuration of type
(A).

\begin{lemma}
Let $p\colon P\to Y$ be a genus $0$ fibration which is smooth outside
$Z_1$ so that the inverse image $p\1(Z_1)$ of $Z_1$ is a simple normal
crossing divisor in $P$. For $s\colon Y\to P$ a section and $C$ a
component of $Z_1$, consider the blowup $P'\to P$ along $s(C)$.

Then $p'\colon P'\to Y$ is again a genus $0$ fibration which is smooth
outside $Z_1$ so that the inverse image ${p'}\1(Z_1)$ of $Z_1$ is a
simple normal crossing divisor, and the birational transform (or strict 
transform) of $s(Y)$ gives a section of $P'$ over $Y$. Moreover, if $t$
is a section of $p$ which is disjoint from $s(C)$ then the birational
transform of $t(Y)$ continues to be a section of $p'$.
\end{lemma}

\begin{proof}
The first two statements and the final statement about $P'$ are obvious.
The union of $s(Y)$ with $p\1(Z_1)$ is a simple normal crossing divisor
since $s$ is a section. Moreover, $s(C)$ is the locus of intersection of
$s(Y)$ and one of the components of $p\1(C)$. Thus the blowup
preserves the property of being a simple normal crossing divisor.
\end{proof}

Suppose that $n_P>1$ and let $N_P$ be the number of $C'$ attaining
this maximum. For each such $C'$ and each pair $i$, $j$ such that
$s_i(C)=s_j(C)=C'$ let $m(i,j,C')$ denote the multiplicity of
intersection of $s_i(Y)$ and $s_j(Y)$ along $C'$. We are in the
situation of the following lemma

\begin{lemma}
Let $P$ be a smooth variety, $D_1$ and $D_2$ smooth divisors meeting
along a smooth codimension $2$ locus $B$ with multiplicity $m>0$. Let
$P'$ be the blowup of $P$ along $B$, and let $D_1'$ and $D_2'$ be the
birational transforms of $D_1$ and $D_2$ respectively. Then the
multiplicity of intersection of $D_1'$ and $D_2'$ is $m-1$ along a
codimension two locus $B'$ lying over $B$.
\end{lemma}

\begin{proof}
In a neighbourhood of the generic point of $B$, the given condition can
be written as follows
 \[
 \Oh_P(D_1)\tensor\Oh_{D_2} = \Oh_{D_2}(m\cdot B)
 \]
Let $E$ denote the exceptional divisor of the blowup $\ep\colon P'\to P$
and let $B'$ be the intersection $D_2'\cap E$; this is a section of the
$\proj^1$-bundle $E\to B$. Now the birational transform of $D_1$
represents the divisor class $\ep^*(D_1)-E$ and thus we see that
 \[
 \Oh_{P'}(D_1')\tensor\Oh_{D_2'} =
 (\Oh_P(D_1)\tensor\Oh_{D_2'}) \tensor \Oh_{D_2'}(-B')
 = \Oh_{D_2'}((m-1)\cdot B'),
 \]
which proves the result.
\end{proof}

 From this lemma, we see that the multiplicity of intersection of the
birational transforms of $s_i(Y)$ and $s_j(Y)$ is $m(i,j,C')-1$. If this
becomes zero then $n(C')$ drops, hence either $N_P$ decreases or $n_P$
does so. This completes the argument by induction since all these
numbers are positive and we wish to obtain the situation where $n_P=1$.

Now, we replace $X$ by the normalisation of the fibre product
$X\times_{\proj_Y(F)} P$ where $P\to Y$ is the configuration of type
(A) with $n_P=1$ obtained above. We then have a finite morphism
$h\colon X\to P$. The image $h(Z)$ of $Z$ and the branch locus of $h$
are contained in the simple normal crossing divisor consisting of the
finitely many disjoint sections of the configuration (A) and
$p\1(Z_1)$. This proves the Claim~\ref{claim2}.

\section{Toric singularities}\label{toric}
As a last step we must prove the result in the following situation.
There is a finite map $f\colon X\to W$ with $W$ smooth and $X$ normal and
$D=\bigcup_{i=1}^n D_i$ a divisor with simple normal crossings so that
$f$ is \'etale outside $D$; moreover $Z$ is a union of (some of) the
components of $f\1(D)$. So we need to construct a birational
morphism $\ep\colon \wave X \to X$ so that $\wave X $ is smooth and
$\ep\1f\1(D)$ is a divisor with simple normal crossings in
$\wave X$.

We will show (see Lemma~\ref{strictness} below) that the inclusion
of $X\setminus f\1(D)$ in $X$ is a strict toroidal embedding in the
sense of \cite{Kempf-etal}, where the desingularisation problem for
such embeddings has been studied and solved. For the sake of
completeness we also give a brief summary of their method. Many
statements below are given without detailed proofs---readers
are invited to complete the arguments on their own or look for proofs in
one of the books on toric geometry such as Oda~\cite{Oda} or
Fulton~\cite{Fulton}.

\subsection{Affine toric singularities}\label{affine}
Let us first consider the simple situation where $W=\A^n$ and $D$ is a
union of coordinate hyperplanes. In this case $W\setminus D$ is
isomorphic to $\G_m^r\times \A^{n-r}$ where $r$ is the number of
components of $D$. Hence (because we are working over the complex
numbers $\C$) its fundamental group is a product of infinite cyclic
groups and the finite cover $X\setminus f\1(D)$ is also isomorphic to
$\G_m^r\times \A^{n-r}$. In fact, the homomorphism $f^*$ on rings has
the form
\begin{eqnarray}
\nonumber
\C[z_1,z_1\1,\dots,z_r,z_r\1,z_{r+1},\dots,z_n] & \to &
  \C[t_1,t_1\1,\dots,t_r,t_r\1,t_{r+1},\dots,t_n] \\
\label{gentor}
 z_i & \mapsto &
 \begin{cases}
 m_i & \text{~for $i\le r$}\\
 t_i & \text{~for $i>r$}
 \end{cases}
 \end{eqnarray}
where the $m_i$ are certain monomials (with negative powers allowed) in
$t_1,\dots,t_r$. The natural action of the torus $\G_m^n$ on $X\setminus
f\1(D)$ then descends to an action via $f$ making this map equivariant.
Moreover, since $X$ is the normalisation of $W$ in $X\setminus f\1(D)$,
the action extends to $X$. Since the map $f$ is equivariant there are
only finitely many orbits for the torus action on $X$; in other words,
$X$ is an (affine) {\em toric variety}.

An explicit description of $X$ can be given as follows. Let $M$ be the
free Abelian group of all monomials in the variables $t_1,\dots,t_r$. Let
$M^+$ be the {\em saturated} submonoid of $M$ generated by the $m_i$. Then
$X=\Spec \C[M^+]\times\A^{n-r}$. For future reference, we note that there
is a unique closed orbit in $X$ which maps isomorphically to the closed
orbit $0\times \A^{n-r}$ in $W$. Moreover, $M$ is the group of Cartier
divisors supported on $f\1(D)$ and $M^+$ is the submonoid of effective
Cartier divisors.

Let $N$ be the dual Abelian group to $M$ and 
 \[
 N^+=\bigl\{ n \,\bigm|\, \text{$n(m)\ge0$ for all $m\in M^+$}\bigr\}
 \]
the ``dual'' monoid to $M^+$. Then $N^+$ is a finitely generated
saturated monoid in $N$ (like $M^+$ in $M$). Let $\si$ be any finitely
generated submonoid of $N^+$ and
 \[
 M^{\si}=\bigl\{ m \,\bigm|\,
 \text{$n(m)\ge0$ for all $n\in N^+$}\bigr\}
 \]
the dual monoid in $M$ (which contains $M^+$). Then
$X(\si)=\Spec\C[M^{\si}]\times\A^{n-r}$ is an affine toric variety and
the natural morphism $X(\si)\to X$ is birational and equivariant for the
torus action. Moreover, we see that $X(\si)$ is non\-singular and the
pullback to $X(\si)$ of $D$ is a simple normal crossing divisor if and
only if the monoid $\si$ is generated over nonnegative integers by a
(sub-)basis of $N$; such a monoid is called {\em simplicial}.

Let $\Si=\{\si_i\}_{i=1}^n$ be a collection of finitely generated
saturated submonoids of $N^+$ which give a subdivision of $N^+$; i.e.,
$N^+$ is the union of all the $\si_i$ and for any pair $\si_i$, $\si_j$
their intersection is a $\si_k$ for some $k$. We then obtain a collection
of equivariant birational morphisms $X_i=X(\si_i)\to X$ so that if
$\si_j\subset\si_i$ then $X_j\subset X_i$ in a natural way. Thus we can
patch together the $X_i$ to obtain $X_{\Si}\to X$ which is birational and
equivariant (but $X_{\Si}$ need not be affine any more). Moreover, the
condition that the $\si_i$ cover $N^+$ implies that
$X_{\Si}\to X$ is proper.

Thus to obtain a desingularisation of $X$ it is enough to find a
subdivision of $N^+$ consisting entirely of simplicial monoids; an easy
enough combinatorial problem solved by barycentric subdivision. The
intrepid reader is warned that proving that the resulting morphism is
projective is a little intricate since an arbitrary simplicial
subdivision {\em need not} result in a projective morphism; however, the
barycentric subdivision does yield a projective morphism.

\subsection{Local toric singularities}\label{local}
Now we examine the general case locally. Let $x\in X$ be any point such
that $w=f(x)$ lies in $D$. There is an analytic neighbourhood $U$ of $w$
in $W$ and coordinates on $U$ so that $D\cap U$ is given by the
vanishing of a product of coordinate functions; by further shrinking
$U$ we can assume that $U$ is a polydisk in these coordinates. Let $V$
be the component of $f\1(U)$ which contains $x$. The normality of $X$
implies the normality of the analytic space $V$. Hence, the open subset
$V\setminus f\1(D)$ is connected. So it is a topological cover of
$U\setminus D$. The coordinate functions give an inclusion $U\into
W'=\A^n$ so that $D\cap U$ is the restriction to $U$ of a union of
coordinate hyperplanes in $W'$. The resulting inclusion $U\setminus
D\into \G_m^r\times\A^{n-r}$ induces an isomorphism of fundamental
groups. Thus there is a Cartesian square in which the horizontal maps are
inclusions:
 \[
\renewcommand{\arraystretch}{1.5}
\begin{matrix}
V\setminus f\1(D) & \longrightarrow & \G_m^r\times \A^{n-r} \\
f \downarrow \hphantom{f} &&  \hphantom{f'} \downarrow  f' \\
U\setminus D  & \longrightarrow & \G_m^r\times \A^{n-r}
\end{matrix}
 \]
and $f'$ is a covering of the form \eqref{gentor} above. Let $X'$
denote the normalisation of $W'$ in the cover $f'$. By the normality
of $X$ (and hence the normality of $V$ as an analytic space) it
follows that we obtain a commutative diagram
 \[
\renewcommand{\arraystretch}{1.5}
\begin{matrix}
V & \longrightarrow & X'\\
f \downarrow \hphantom{f} &&  \hphantom{f'} \downarrow  f' \\
U & \longrightarrow & W'
\end{matrix}
 \]
so that $V$ is isomorphic to an analytic open neighbourhood of a point
in the toric variety $X'$. By means of the two diagrams above we can
carry over the desingularisation of \ref{affine} to the local analytic
space $V$.

\subsection{General toroidal embeddings}
The local description \ref{local} can be repeated in a suitable
neighbourhood of any point $x\in X$. This shows that the inclusion of
$X\setminus f\1(D)$ in $X$ is a {\em toroidal embedding}. The
desingularisations obtained locally need to be constructed in a coherent
manner so that they ``patch up''.

Let us stratify $W$ by connected components of intersections of the form
\begin{equation}\label{intn}
 (D_{i_1}\cap\dots\cap D_{i_r})\setminus \bigcup_{j\neq i_s} D_j.
\end{equation}
If $U$ and $V$ are as in \ref{local}, then there is a unique stratum
$S$ of $W$ so that $S\cap U$ is closed. Under the inclusions of $U$ in
$W'$ and $V$ in $X'$ of \ref{local} the strata correspond to the
orbits of the torus action. As we have noted in \ref{affine} the
unique closed strata in $U$ and $V$ then become isomorphic. Thus, if
$T=f\1(S)$ then $T\cap V$ is closed in $V$ and the morphism $T\cap V\to
S\cap U$ is an isomorphism; indeed $T\cap V$ and $S\cap U$ are the
restrictions of the closed orbits in $X'$ and $W'$ respectively. Thus we
see that if $S$ is {\em any} stratum in $W$ and $T$ a connected component
of $f\1(S)$ then $T\to S$ is \'etale and proper. We thus stratify $X$ by
connected components of the inverse images of strata in $W$ under $f$.
Let $\{T_a\}_{a\in A}$ denote this stratification; by abuse of notation
we define $S_a=f(T_a)$.

Let $E=\bigcup_{j=1}^m E_j$ be the decomposition into irreducible
components of the inverse image $f\1(D)$ of $D$. Then there is a
function $i\colon \{1,\dots,m\}\to\{1,\dots,n\}$ so that
$f(E_j)=D_{i(j)}$ for all $j$. For any stratum $T_a$ let $X^a$ denote
the complement in $X$ of all those $E_j$ that do not meet $T_a$.
Similarly, we denote by $W^a$ the complement in $W$ of all $D_i$ that do
not meet $S_a$. The morphism $f$ clearly maps $X^a$ to $W^a$.

\begin{lemma}\label{strictness}
Let $M_a$ be the Abelian group of all Cartier divisors in $X^a$ with
support in $E\cap X^a$; let $M^+_a$ be the submonoid consisting of
effective Cartier divisors. Then $M_a$ has rank $r_a=\codim_X T_a$. The
distinct analytic components of $E$ in a neighbourhood of any point $x$
of $X$ are precisely the algebraic components; i.e., $X\setminus E\into
X$ is a {\em strict} toroidal embedding.
\end{lemma}

\begin{proof}
Let $M_a'$ be the group of Cartier divisors on $W^a$ with support on
$D\cap W^a$. Clearly, this is the free Abelian group on those $D_i$ which
contain $S_a$. Since $S_a$ is a connected component of an intersection
of the type \eqref{intn}, there are exactly $r_a=\codim_W S_a$ of such
$D_i$. The pullback homomorphism makes $M'_a$ a subgroup of $M_a$. On the
other hand, consider an analytic neighbourhood $V=f\1(U)$ of some point
$x\in T_a$ as in \ref{local}. If $M_x''$ denotes the group of all
Cartier divisors on $V$ supported on $E\cap V$; then we have noted in
\ref{affine} that $M_x''$ is a free Abelian group of rank $r_a$.
Moreover, $M_a$ is included as a subgroup of $M_x''$ under the
restriction from $X_a$ to $V$.

It follows that the homomorphism $M_a\to M_x''$ has finite cokernel; but
then the normality of $X$ means that it is surjective. In particular, we
see that the distinct analytic components of $E\cap V$ correspond to
distinct algebraic components of $E$. This concludes the proof.
\end{proof}

Let $M_a^+$ denote the monoid of effective divisors in $M_a$; under the
isomorphism $M_a\to M_x''=M$ this maps isomorphically onto the submonoid
$M^+$ considered in \ref{affine}. Thus $X_a$ is smooth (and the
divisor $E$ is a simple normal crossing divisor) if and only if $M_a^+$
is simplicial. If $T_b$ lies in the closure of $T_a$ then $X_a$ is an
open subset of $X_b$. This induces by restriction a (surjective)
homomorphism $M_b\to M_a$ which further restricts to a surjection
$M_b^+\to M_a^+$. Thus we have a (finite) projective system of monoids.

Let $N_a$ and $N_a^+$ be the dual objects as defined in \ref{affine}.
These form a finite injective system of monoids. By a compatible family
$\sS$ of subdivisions we mean a subdivision $\Si_a$ of $N_a^+$ for each
$a$ so that the subdivision $\Si_b$ restricts to the subdivision $\Si_a$
on the submonoid $N_a^+$ of $N_b^+$. We then obtain a proper birational
morphism $X_{\Si_b}\to X_b$ for each $b$ which restricts to $X_{\Si_a}\to
X_a$ on the open subset $X_a$ of $X_b$. Thus, we see that any such
compatible family leads to a proper birational morphism $X_{\sS}\to X$.

Thus, in order to desingularise, we have to find a compatible family of
subdivisions so that each of the new monoids is simplicial. This is
achieved by the barycentric subdivision. As seen earlier, this ensures
that the morphism $X_{\sS}\to X$ is locally projective. Since $X$ is
projective, this morphism is indeed projective.

\nocite{Hironaka}


\begin{thebibliography}{1}

\bibitem{Ab-Jong}
D.~Abramovich and A.~J.~de~Jong, {\em Smoothness, semistability and
 toroidal geometry},
 J. Algebraic Geom. {\bf 6} (1997), no. 4, 789--801;
 {\tt arXiv:math.AG/9603018}, 1996.

\bibitem{Belyi}
G.~V.~Belyi, {\em Galois extensions of a maximal cyclotomic field},
Izv.~Akad.~Nauk SSSR Ser.~Math.\ {\bf 43} (1979), no. 2, 267--276.

\bibitem{Bog-Pant}
F.~Bogomolov and T.~Pantev, {\em Weak Hironaka theorem},
 Math. Res. Lett. {\bf 3} (1996), no. 3, 299--307; 
{\tt arXiv:math.AG/9603019}, 1996.

\bibitem{Fulton}
W.~Fulton, {\em Introduction to Toric varieties}, Annals of Math
Studies, Vol. 151, Princeton University Press, Princeton, USA, 1993.

\bibitem{Hironaka}
H.~Hironaka, {\em Resolution of singularities of an algebraic variety over a
 field of characteristic zero: {I}, {II}}, Ann.\ of Math.(2) {\bf 79}
(1964), 109--326.

\bibitem{Jong}
A.~J.~de~Jong, {\em Smoothness, semistability and alterations},
 Publ. Math. IHES {\bf 83} (1996), 51--93.

\bibitem{Kempf-etal}
G.~Kempf, F.~Knudsen, D.~Mumford, and B.~{Saint-Donat}, Lecture Notes in
Mathematics, no. 339, Springer-{V}erlag, Berlin-Heidelberg-New York,
1973.

\bibitem{Redbook}
D.~Mumford, {\em The Red book of varieties and schemes}, Narosa
Publishers, New Delhi, India, 1996, Indian edition.

\bibitem{Oda}
T.~Oda, {\em Lectures on Torus Embeddings and applications}, TIFR
Lecture Notes, TIFR, Bombay, India, 1982.

\end{thebibliography}
\end{document}